\documentclass[3p,times]{article}
\usepackage{amsmath}
\usepackage{amsfonts}
\usepackage{mathrsfs}
\usepackage{amssymb}
\usepackage{times}
\usepackage{subcaption}
\usepackage{xcolor}
\usepackage{authblk}
\usepackage{graphicx}%
\newtheorem{theorem}{Theorem}
\newtheorem{proposition}[theorem]{Proposition}
\newtheorem{remark}[theorem]{Remark}

\newtheorem{corollary}[theorem]{Corollary}
\newtheorem{lemma}[theorem]{Lemma}

\title{Dependence and mixing for perturbations of copula-based Markov chains}
\author[1]{Martial Longla\thanks{mlongla@olemiss.edu}}
\author[1]{Mathias Muia Nthiani \thanks{mnmuia@go.olemiss.edu}}
\author[2]{Fidel Djongreba Ndikwa.\thanks{fideldjong@gmail.com}}

\affil[1]{University of Mississippi, Department of mathematics}
\affil[2]{University of Maroua, Department of mathematics}
\counterwithin{equation}{section}
\counterwithin{theorem}{subsection}
\begin{document}
\maketitle

\begin{abstract}
This paper explores the impact of perturbations of copulas on dependence properties of the Markov chains they generate. We use an observation that is valid for convex combinations of copulas to establish sufficient conditions for the mixing coefficients $\rho_n$, $\alpha_n$ and some other measures of association. New copula families are derived based on perturbations of copulas and their multivariate analogs for $n$-copulas are provided in general. Several families of copulas can be constructed from the provided framework.  
\end{abstract}


\bigskip

\section{Introduction}
Many authors have contributed to constructing new families of copulas to model dependence among variables or factors in economics, finance, risk management and other applied fields. Following the ideas of Durante and al. (2013), we consider the perturbation method that adds to a copula an extra term called perturbation. This modification is often used to increase flexibility of the dependence structure of the considered models. Here, we look at other classes of modifications and their impact on the dependence structure as studied by Komornik and al. (2017). We consider the long run impact of these perturbations on the dependence structure and the measures of association. The case is presented for $\rho$-mixing and $\alpha$-mixing. This is a continuation of the work done by Longla and al. (2021) in which $\psi$-mixing, $\beta$-mixing and $\phi$-mixing were explored. The paper is structured as follows. Introduction in Section 1 is divided into several parts, each of which is dedicated to a specific topic of interest. Copulas are introduced in Section 1.1, perturbations of copulas are defined in Section 1.2 with some results on perturbations, Section 1.3 is dedicated to the perturbation of components of a vector and their impact on the copulas,  Section 1.4 is dedicated to the mixing coefficients of interest, Section 1.5 visits the impact of perturbations on some dependence coefficients (measures of association). In Section 2, we present the impact of perturbations on the mixing structure of copula-based Markov chains, addressing $\alpha-$mixing in Section 2.1 and $\rho$-mixing in Section 2.2. In Section 3, we provide the proofs of the main results of this work.

\subsection{Copula theory}
The importance of copulas is due to the fact that they are tools that link univariate distributions to the joint distribution of the random variables. They were first presented in Sklar (1959), and over the recent years they have become very important for analyzing temporal dependence of time series. For $n\geq 2$, an $n$-copula is an $n$-variate distribution function on $[0,1]^{n}$ with uniform marginals. One of the most used types of copula is the $2$-copula, that is mostly used for copula-based Markov chains. By definition, a $2$-copula is a bivariate distribution function that has uniform marginal distributions on $I=[0,1]$. They have gained a lot of popularity in financial, economics and risk management time series because they are used to generate Markov chains with important needed properties. For example, Ibragimov (2009) illustrates applications of copulas in Econometric theory. The definition of a $2$-copula and related topics can be found in Nelsen (2006). Formally, a function $C: [0,1]^{2}\rightarrow [0,1]$ is called bivariate copula if it satisfies the following conditions:
\begin{enumerate}
	\item[i.] $C(0,x)=C(x,0)=0$ (meaning that $C$ is grounded); \item[ii.]$C(x,1)=C(1,x)=x, \forall x\in I$ (meaning each coordinate is uniform on $I$);
\item[iii.] $C(a,c)+C(b,d)-C(a,d)-C(b,c)\geq 0, \forall\
[a,b]\times[c,d]\subset I^{2}.	$
\end{enumerate}
Darsaw and al. (1992) derived the transition probabilities for stationary Markov chains with uniform marginals on $[0,1]$ as $P(X_{n}\in A|X_{n-1}=x)=C_{,1}(x,y), \forall n\in\mathbb{Z}$ and $A=\left(-\infty,y\right]$. Here, $C_{,i}(x,y)$ denotes the derivative with respect to the $i^{th}$ variable. This property has been used by many authors to establish mixing properties of copula-based Markov chains. We can cite Longla (2015), Longla (2014), Longla and Peligrad (2012) who provided some results for reversible Markov chains and Beare (2010) who presented results for $\rho$-mixing.  We are assuming in the sequel that variables in a copula-based Markov chain have the uniform distribution on $[0,1]$.

It's been shown in the literature (see Darsow and al. (1992) and the references therein) that if $(X_1, \cdots, X_n)$ is a Markov chain with consecutive copulas $(C_1, \cdots, C_{n-1})$, then for any integer $k$, the fold product given by 
$$C(x,y)=C_k*C_{k+1} (x, y)=\int^1_0 C_{k,2}(x, t)C_{k+1,1}(t, y)dt$$ is the copula of $(X_k,X_{k+2})$ and the $\star$-product given by
$$ C(x,y,z)=C_k\star C_{k+1} (x, y,z)=\int_0^y C_{k,2}(x, t)C_{k+1,1}(t, z)dt$$ is the copula of $(X_k,X_{k+1},X_{k+2})$. The $n$-fold product of $C(x,y)$ denoted $C^n(x,y)$ is defined by the recurrence $C^{1}(x,y)=C(x,y)$, $C^{n}(x,y)=C^{n-1}*C(x,y), \quad \text{for $n>1.$}$ This product is the joint distribution of $(X_k, X_{k+n})$ when all copulas are equal along the chain. 
The most popular copulas are $\Pi(u,v)=uv$  (the independent copula), the Hoeffding lower and Upper bounds  $W(u,v)=\max(u+v-1,0)$ and $M(u,v)=\min(u,v)$ respectively.  Convex combinations of copulas $\{C_1(x,y), \cdots, C_k(x,y)\}$ defined by $\displaystyle \{ C(x,y)=\sum_{j=1}^{k}a_j C_j(x,y), 0\leq a_j, \sum_{j=1}^{k} a_j=1\}$ are also copulas and the following holds.

\begin{proposition}{convex combinations and fold product} \label{convex}
\begin{enumerate}
\item The fold product of convex combinations of copulas is a convex combination of copulas.
\item The $n$-fold product of any convex combination of copulas is a convex combination of copulas.
\end{enumerate}
\end{proposition}
The proof of this proposition relies on the fact that the fold product defines a bilinear operator. For any copulas $A_1(x,y), A_2(x,y), B_1(x,y), B_2(x,y)$ be copulas. It is easy to verify that 
$$A_1*(a_1 B_1+a_2 B_2)(x,y)=a_1 A_1*B_1 (x,y)+a_2 A_1*B_2(x,y)$$
$$(a_1 B_1+a_2 B_2)*A_1(x,y)=a_1 B_1*A_1 (x,y)+a_2 B_2*A_1(x,y)$$
So, $(a_1 A_1+a_2 A_2)*(b_1 B_1+b_2 B_2)(x,y)=a_1b_1A_1*B_1 (x,y) + a_1b_2A_1*B_2 (x,y) + a_2b_1A_2*B_1 (x,y) + a_2b_2A_2*B_2 (x,y).$
if all initial coefficients are positive with $a_1+a_2=1$ and $b_1+b_2=1$, then the last set of coefficients is made of positive numbers with $a_1b_1+a_1b_2+a_2b_1+a_2b_2=a_1(b_1+b_2)+a_2(b_1+b_2)=1$. Thus, the fold product is a convex combination. The second part follows by multiplying the convex combination by itself. This result by itself is elementary, but plays a big role in establishing many results or observations in the sequel.
From the proof of Proposition \ref{convex}, it follows that for a convex combination $\displaystyle C(x,y)=\sum_{i=1}^{k} a_i C_i(x,y)$ and any integer $n$ the following formula holds. 
\begin{equation} \label{convobs}
C^{n}(x,y)=\sum_{j=1}^{k^n} {b_j}\times  {_{1}C_{j}}*\cdots *{_{n}C_{j}(x,y)}, 
\end{equation}
where $\displaystyle \sum_{j=1}^{k^n}b_j=1, \quad b_j\ge 0$ , each of the copulas $_{i}C_{j}(x,y)=C_{j_i}(x,y)$ for some $j_i\in\{1, \cdots, k\}$ and the sum is over all possible products of $n$ copulas selected from the original $k$ copulas with replacement. The notation $_{i}C_{j}$ indicates that the copula $C_{j_i}$ was selected in the given $j^{th}$ element of $B=\{ C_1, \cdots, C_k\}^{n}$. 
If we denote $c^n(x,y)$ the density of the copula $C^{n}(x,y)$, then the following holds.
 
\begin{lemma} \label{bound}
If a copula $C(x,y)$ is such that its density $c(x,y)$ satisfies $c(x,y)\ge K$ on a set of Lebesgue measure 1 for some real number $K>0$, then for all positive integer $n$ we have $c^{n}(x,y)\ge K^n$ and $C^{n}(x,y)\ge xyK^n$. 
\end{lemma}

\subsection{Perturbations of copulas}
Copulas are well known tools for understanding dependence among random variables. They offer a choice of appropriate models for dependence between random variables independently from selection of marginal distributions. Often in model analysis it is necessary to understand the dependence structure of random variables for inference purposes. Commonly used measures of dependence include Kendall's tau ($\tau$), Spearman's rho ($\rho$), Blomqvist's beta ($\beta$), Gini's gamma ($\gamma$), among others. Sometimes linear correlation fails to be reliable in measuring dependence between random variables; especially where non-linear monotonic transformations of random variables are involved. It is at this point that the concordance and dependence measures (e.g. Kendal's $\tau$, Spearman's $\rho$) perform best as they are invariant under monotonic transformations.  We consider the copula $C$ and the related copula $C_H$ defined by $C_H(u,v)=C(u,v)+H(u,v),$
where $H : [0, 1]^2 \to \mathbb{R}$ is a continuous function. The function $H$ is the perturbation factor and the copula $C_H$ is called a perturbation of $C$.  This perturbation method introduced by Durante and al. (2013) and its special case with $$H(x,y)=H_\theta(u,v):=\theta (u-C(u,v))(v-C(u,v)), \theta\in [0,1]$$ was investigated by Mesiar and al. (2015). We  consider the class of perturbations defined in general by 
\begin{equation}
C_{\theta}(u,v)=C(u,v)+ \theta[C_1(u,v)- C(u,v)], \quad \theta\in [0,1], \quad C_1 \mbox{ is a copula,}
\end{equation}
then classes of perturbations with $C_1(u,v)=\Pi(u,v)=uv$ and $C_1(u,v)=M(u,v)$ are respectively as follows
\begin{equation}
C_{\theta,\Pi}(u,v)=C(u,v)+ \theta[ uv- C(u,v)]=(1-\theta)C(u,v)+\theta \Pi(u,v), \quad \theta\in [0,1],  \label{pi}
\end{equation}
\begin{equation}
C_{\theta,M}(u,v)=C(u,v)+ \theta[ M(u,v)- C(u,v)]=(1-\theta)C(u,v)+\theta M(u,v), \quad \theta\in [0,1].  \label{M}
\end{equation}

This class of perturbations modifies all copulas except when $C=C_1$.  Technically, perturbations of this kind happen when the bivariate population of interest with copula $C(u,v)$ has been contaminated by another population with copula $C_{1}(u,v)$. We consider the long run effect of the perturbations given by \eqref{pi} and \eqref{M}. By Longla and al. (2021), we have 
\begin{equation}
C^{n}_{\theta, M}(u,v)= \sum_{i=1}^{n}{n \choose i}\theta^{n-i} (1-\theta)^{i}C^{i}(u,v)+\theta^nM(u,v) \label{Mcop}
\end{equation} 
\begin{equation}
C^{n}_{\theta, \Pi}(u,v)=(1-\theta)^nC^{n}(u,v)+(1-(1-\theta)^n)\Pi(u,v). \label{Pcop}
\end{equation} 
Based formulas \eqref{Mcop} and \eqref{Pcop}, via simple computations, the following is established.
\begin{proposition} \label{MP}
For any $\theta\in(0,1)$, 
\begin{enumerate} 
\item If $\displaystyle \lim_{n\to\infty} C^n(u,v)=C_o(u,v)$, then
$\displaystyle
\lim_{n\to\infty} C^n_{\theta, M}(u,v)=C_o(u,v)
$
\item $\displaystyle \lim_{n\to\infty} C^n_{\theta, \Pi}(u,v)=\Pi(u,v). $
\end{enumerate}
\end{proposition}

The proof of Proposition \ref{MP}  relies on  Proposition 1.0.1 of Longla and al. (2021). It makes more sense looking at this proposition to consider formula \eqref{Pcop} as a perturbation of the independence copula by means of the copula $C(u,v)$. As an example, we can take a copula $C(u,v)$ from the two parameters Mardia family of copulas (a subfamily of this class of copulas is the one parameter Frechet copula family ). Copulas from this family are defined as follows.
\begin{equation}\label{Mardia}
C=aM(u,v)+bW+(1-a-b)\Pi, \quad \mbox{for some}\quad a,b,a+b\in[0,1].
\end{equation}
When $a=\theta^2(1+\theta)/2$ and $b=\theta^2(1-\theta)/2$ the Mardia copula turns into a one parameter copula called Frechet copula.
Longla (2014) has shown that for this class of copulas, $$C^n=\frac{1}{2}[(a+b)^n+(a-b)^n]M+\frac{1}{2}[(a+b)^n-(a-b)^n]W+[1-(a+b)^n]\Pi.$$
Therefore,
\begin{equation}
\lim_{n\to\infty}C^n(u,v)=\begin{cases} \frac{1}{2}M(u,v)+\frac{1}{2}W(u,v), & \mbox{if} \quad a+b=1, a\ne1,b\ne 1\\ \mbox{Undefined}, & \mbox{if} \quad b=1, a=0, \\ M(u,v), & \mbox{if}\quad a=1, b=0, \\ \Pi(u,v), & \mbox{if} \quad a+b<1. \end{cases}
\end{equation}
Convergence in this case is uniform convergence. This can be obtained via simple computations that we omit here. This limit by itself shows that copulas from the Mardia family generate Markov chains for which $(X_0,X_n)$ are either asymptotically independent (when $a+b<1$), or asymptotically almost surely equal (when $a=1, b=0$) or asymptotically almost surely $X_n=X_0$ with probability $1/2$ or $X_n=1-X_0$ with probability $1/2$ (when $a+b=1, a\ne1, b\ne 1$). For $b=1, a=0$ there is no limit of the copula sequence. This limit also shows that Frechet copulas either converge to $\frac{1}{2}M(u,v)+\frac{1}{2}W(u,v)$ (for $|\theta|<1$) or converge to $M(u,v)$ (for $\theta=1$) or have no limit when $\theta=-1$. The last case is because $W^2(u,v)=M(u,v)$, and $M*W(u,v)=M*W(u,v)=W(u,v)$.

As mentioned by Durante and al. (2013) perturbations of copulas have been studied under simplified assumptions on $C$ and $H$. In particular, in the literature it is often assumed that $C(u,v)=\Pi(u,v).$ The model is interpreted as a contamination of independence. 
We will also consider in this paper another type of perturbation that was studied in Longla and al. (2021) and was considered for measures of association in the works of Durante and al (2013),  Mesiar  and al (2015) and Komornik and al (2017). This perturbation method acts at the level of the random variables, polluting the underlying variable $(X)$ with some noise variable $(Z)$. In the cited works, the variable $X$ is mostly bivariate and the variable Z is also bivariate, independent of $X$ and with possible dependence between its components. In this work, we will consider a vector $X=(X_1, \cdots, X_k)$ and a perturbation vector $Z=(Z_1, \cdots, Z_k)$ for some positive integer $k$. The goal is to derive the copula of $X+Z$ under various conditions on $X$ and $Z$ and provide some examples for $X$ having copula or distribution $M(u_1, \cdots, u_k)$ or $\Pi(u_1, \cdots, u_k)$.
\subsection{Perturbation of $X$ via $X+Z$ and its copulas}
Assume we have a random $k$-dimensional vector $X=(X_1, \cdots, X_k)$ of copula $C(u)$ and marginal distributions $F_1, \cdots, F_k$. Consider a noise vector $Z=(Z_1, \cdots, Z_k)$ independent of $X$. We are interested in extracting the copula of $X+Z$ from its joint distribution function.  Assume $\tilde{F}_k$ is the cumulative distribution function of the random variable $X_k+Z_k$ and $u=(u_1, \cdots, u_k)$. 

\begin{theorem}\label{IndXZpart}
If $Z=(Z_1, \cdots, Z_s, 0\cdots, 0)$ for $s\le k$ and $(Z_1, \cdots, Z_s)$ has independent components with marginal distributions $G_1, \cdots, G_s$, then the copula of $X+Z$ is 
\begin{equation} 
C_1(u)=\int_{-\infty}^{ \infty}\cdots\int_{-\infty}^{ \infty} H(u,t)dG_1(t_1)\cdots dG_s(t_s), 
\end{equation}
where $H(u,t)=C(F_1(\tilde{F}_1^{-1}(u_1)-t_1), \cdots, F_s(\tilde{F}_s^{-1}(u_s)-t_s), u_{s+1}, \cdots, u_k)$.
\end{theorem}
A consequence of Theorem \ref{IndXZpart} is the following:

\begin{corollary}\label{IndXZ}
If $Z=(Z_1, \cdots, Z_k)$ has independent components with marginal distributions $G_1, \cdots, G_k$, then the copula of $X+Z$ is 
\begin{equation} 
C_2(u)=\int_{-\infty}^{ \infty}\cdots\int_{-\infty}^{ \infty} C(F_1(\tilde{F}_1^{-1}(u_1)-t_1), \cdots, F_k(\tilde{F}_k^{-1}(u_k)-t_k))dG_1(t_1)\cdots dG_k(t_k).
\end{equation}
\end{corollary}

\begin{theorem}\label{XZpart}
If $Z=(W, \cdots, W, 0\cdots 0)$ where $W$ has distribution $G$ and $Z_i=0$ for all $i>s$, then the copula of $X+Z$ is 
\begin{equation} 
C_3(u)=\int_{-\infty}^{ \infty} C(F_1(\tilde{F}_1^{-1}(u_1)-t), \cdots, F_s(\tilde{F}_s^{-1}(u_s)-t), u_{s+1},\cdots, u_k)dG(t).
\end{equation}
\end{theorem}

A consequence of Theorem \ref{XZpart} is the following 
\begin{corollary} \label{XZ}
If $Z=(Z_1, \cdots, Z_1)$ where $Z_1$ has distribution $G$, then the copula of $X+Z$ is 
\begin{equation}
C_4(u)=\int_{-\infty}^{ \infty} C(F_1(\tilde{F}_1^{-1}(u_1)-t), \cdots, F_k(\tilde{F}_k^{-1}(u_k)-t))dG(t).
\end{equation}
\end{corollary}
The proofs of these statements are based on simple computations using definitions of probabilities via conditioning an Sklar's theorem as in Longla and al. (2021).
\begin{remark}
Notice that in each of the two theorems, the position of the variables that are modified doesn't change the details of the proof. In fact, if $\sigma: \{1, \cdots, n\}\to \{1, \cdots, n\}$ is a permutation such that $\sigma(u)=v$, then the copula of $X_\sigma=X+Z_\sigma$, where components of $Z$ have been permuted via $\sigma$ to obtain $Z_\sigma$ is $$C_{i,\sigma}(u)=C_{i}(v), \quad i=1 ,3.$$
\end{remark}

\subsection{Mixing coefficients} 
The concept of mixing is important in the study of statistical properties of deterministic dynamical systems (ergodic theory). The most important mixing coefficients include but are not limited to $\phi$-mixing, $\rho$-mixing, $ \beta$-mixing, $\alpha$-mixing, $\rho*$-mixing, $\psi'$-mixing and $\psi$-mixing. The definitions of these can be found in Bradley (2007). Other works on mixing include but are not limited to Miller (1994), Miller (1995) Gaposhkin (1991) and the references therein. Bradley (1997) presents a proof that every $\psi'$-mixing Markov chain is $\rho^{*}$-mixing. He also includes a detailed list of relationships between different mixing coefficients.  Longla (2015) showed that Markov chains generated by copulas are {\it lower $\psi$-mixing} when the density of the absolutely continuous part of the generating copula is bounded away from zero on a set of Lebesgue measure 1.  Longla and Peligrad (2012) defined the mixing coefficients for an absolutely continuous copula and an absolutely continuous invariant distribution for the states of a stationary Markov chain that they generate. They provided a result on absolute regularity for mixtures of copulas. Archimedean copulas with non-strict generators have been studied by Longla (2014). A condition was established for exponential $\rho$-mixing of copula-based Markov chains in the case  of non-strict Archimedean copulas. A non-strict Archimedean copula is defined by some convex function $\varphi:[0,1]\to [0,1]$ such that $\varphi^{-1}(u)=0 \quad\forall u>1$ via $ C(u,v)=\varphi^{-1}(\varphi(u)+\varphi(u)).$ The $\alpha$-mixing condition is also known as strong mixing condition in the literature (see Bradley (2007) and Bradley (2005)). We borrow the following historical background from Bradley (2005): The strong mixing condition was introduced by Rosenblatt (1956). The $\psi'$-mixing condition was introduced by Bradley (1983), following the work of Blum and al. (1963) on a related coefficient ($\psi$-mixing, not considered here).  The $\rho$-mixing condition was introduced by Kolmogorov and Rozanov (1960).  These coefficients are defined in the literature as follows (see Bradley (2007))
$$\alpha(\mathscr{A},\mathscr{B})=\sup_{B\in \mathscr{B}, A\in \mathscr{A}}|P(A\cap B)-P(A)P(B)|,$$
$$\rho(\mathscr{A},\mathscr{B})=\sup_{f\in L^{2}(\mathscr{A}),g\in L^{2}(\mathscr{B})}corr(f,g),\quad \quad \psi'(\mathscr{A},\mathscr{B})=\inf_{B\in \mathscr{B}, A\in \mathscr{A}}\frac{P(B\cap A)}{P(A)P(B)},$$
where $\mathscr{A}=\sigma(X_{i}, i\leq 0)$, $\mathscr{B}=\sigma(X_{i}, i\geq n)$ and $P$ is the defined probability measure. For Markov chains generated by an absolutely continuous copula (see Longla (2013) or Longla (2015)) these coefficients are 
$$\psi'_{n}:=\psi'(\mathscr{A},\mathscr{B})=\inf_{A,B\in \mathscr{B}}\frac{|\int_{A}\int_{B}c_{n}(x,y)dydx|}{\lambda(A)\lambda(B)},$$
$$\rho_{n}=\sup_{f,g}\{\int^{1}_{0}\int^{1}_{0}c_{n}(x,y)f(x)g(y)dxdy : ||g||_{2}=||f||_{2}=1, \mathbb{E}(f)=\mathbb{E}(g)=0\},$$
$$\alpha_n:=\alpha(\mathscr{A},\mathscr{B})=\sup_{A,B\in \mathscr{B}}|\int_A\int_{B}(c_{n}(x,y)-1)dydx|,$$
where $c_n$ is the copula density of the random variable $(X_{0}, X_{n})$, $\mathscr{B}$ is the Borel sigma field and $\lambda$ is the Lebesgue measure on $(0,1)$. A stochastic process is said to be an $\alpha$-mixing, if $\alpha_{n}\rightarrow 0$; $\psi'$-mixing if $\psi'_{n}\rightarrow 1$ and $\rho$-mixing if $\rho_{n}\rightarrow 0$. The process is exponentially mixing, if the convergence rate is exponential.

\subsection{Perturbations of copulas and dependence coefficients}
 If $C$ is the copula of $(X,Y)$, then the Spearman's rho is $\displaystyle \rho_S(C) = 12 \int_{0}^1\int_{0}^1 C(u,v)dudv-3;$ Kendall's $\tau$ is $\displaystyle \tau(C)=1-4\int_{0}^1\int_{0}^1 C_{,1}(u,v)C_{,2}(u,v)dudv;$ Blomqvist's $\beta$ is $\displaystyle \beta(C)=4C(1/2, 1/2)-1$ and is called the medial correlation coefficient. The Gini's $\gamma$ is $\displaystyle \gamma(C)=4\int_{0}^{1}C(u,u)+C(u,1-u)du;$ the upper tail dependence coefficient is  $\displaystyle \lambda_U(C)=\lim_{u\to 1}\frac{1-2u+C(u, u)}{1-u}$ and the lower tail dependence is $\displaystyle \lambda_L(C)=\lim_{u\to 0}\frac{C(u, u)}{u}$ (see Nelsen (2006) for more on these coefficients). From formulas \eqref{Mcop} and \eqref{Pcop}, it follows:
 
\begin{theorem}\label{asso}
For any copula-based Markov chain generated by a copula $C(u,v)$, 
\begin{enumerate}
\item  $\displaystyle \rho_S(C^n_{\theta,M})=\sum_{i=1}^{n}{n \choose i}\theta^{n-i} (1-\theta)^{i}\rho_S(C^i)+\theta^n$ and $\rho_S(C^n_{\theta,\Pi})=(1-\theta)^n\rho_S(C^n)$.
\item $\displaystyle \beta(C^n_{\theta,M})=\sum_{i=1}^{n}{n \choose i}\theta^{n-i} (1-\theta)^{i}\beta(C^i)+\theta^n$  and  $\beta(C^n_{\theta,\Pi})=(1-\theta)^n\beta(C^n)$. 
\item $ \displaystyle \gamma(C^n_{\theta,M})=\sum_{i=1}^{n}{n \choose i}\theta^{n-i} (1-\theta)^{i}\gamma(C^i)+3\theta^n$ and $\gamma(C^n_{\theta,\Pi})=(1-\theta)^n\gamma(C^n)+2(1-(1-\theta)^n)$.
\item $\displaystyle \lambda_L(C^n_{\theta,M})=\sum_{i=1}^{n}{n \choose i}\theta^{n-i} (1-\theta)^{i}\lambda_L(C^i)+\theta^n$ and $\lambda_L(C^n_{\theta,\Pi})=(1-\theta)^n\lambda_L(C^n)$.
\item $\displaystyle \lambda_U(C^n_{\theta,M})=\sum_{i=1}^{n}{n \choose i}\theta^{n-i} (1-\theta)^{i}\lambda_U(C^i)+\theta^n$ and $\displaystyle \lambda_U(C^n_{\theta,\Pi})=(1-\theta)^n\lambda_U(C^n)$.
\end{enumerate}
\end{theorem}
If we denote $\xi(C)$ any of the considered measures of association other than $\lambda_L$ or $\lambda_U$, then the following holds when $\theta\ne 1$. 
\begin{remark}
Theorem \ref{asso} shows that as $n\to \infty$, $\xi(C^n_{\theta,\Pi})\to \xi(\Pi)$. When the independence copula is perturbed by any copula, the measure of association is asymptotically unchanged. 
 It can also be seen that $\xi(C^n_{\theta,M})$, as $n\to\infty$, depends solely on properties of $C(u,v)$. Using Longla (2021), we can conclude that the limit is equal to $\lim \xi(C^{n})$ when the later exists.
\end{remark}

\section{Perturbations of copulas and mixing coefficients}

\subsection{On $\alpha$-mixing for copula-based Markov chains}

 The convergence rate of this coefficient is related to the central limit theorem via several conditions in the literature (see Merlevede and Peligrad (2020) and the references therein).  From Lemma 3 of Longla and Peligrad (2012), the following holds.
\begin{corollary} \label{al}
Any convex combination of copulas, each of which generates $\alpha$-mixing Markov chains, generates $\alpha$-mixing Markov chains.
\end{corollary}

Following Corollary \ref{al}, after deriving the $n$-th transition kernel for copula-based Markov chains generated by the perturbation copula, we can establish the following.

\begin{theorem} \label{alpha}
When a copula is perturbed using the independence copula $\Pi(u,v)$, the $\alpha$-mixing rate of convergence is multiplied by an exponential factor. For a perturbation that uses the copula $M(u,v)$, $\alpha_n({C_{\theta,M}})\to 0$ when the initial copula generates $\alpha$-mixing Markov chains. In a case when $\lim \alpha_n(C)=c\ne 0$, $C_{\theta, M}$ doesn't generate $\alpha$-mixing.
\end{theorem}

\begin{remark}
As strange as it might sound, this theorem shows that a perturbation of a non $\alpha$-mixing generating copula by means of an $\alpha$-mixing generating copula turns it into an $\alpha$-mixing generating copula. In this specific case, the limit of $C_{\theta,M}^n(u,v)$ is  actually not influenced by $M(u,v)$ for any $\theta\in(0,1)$ and $C_{\theta,\Pi}^n(u,v)$ converges to $\Pi(u,v)$. This also means that in the case of the perturbation of the independence copula, long run realizations of the Markov chain generated by the perturbation copula are independent of the initial state. 
\end{remark}

The following result for ergodic sequences is Theorem 3.4. of Bradley (2005).
\begin{theorem}\label{Brad1}
Suppose $(X_n,n\in Z)$ is a strictly stationary Markov chain which is ergodic and aperiodic.
 If $\alpha_n <1/4$ for  some $n\ge 1$,  then  $\alpha_n\to 0$ (but  not  necessarily exponentially fast) as $n\to\infty$.
\end{theorem}
Based on Theorem \ref{Brad1}, the following holds for convex combinations of copulas.
\begin{theorem}\label{alpha1}
For any set of copulas $C_1(u,v)\cdots C_k (u,v)$, if there exists a subset of copulas $C_{k_1}\cdots C_{k_s}, $ $s\leq k\in \mathbb{N}$ such that $\alpha(\hat{C})<1/4\quad \text{for}\quad \hat{C}=C_{k_1}*\cdots*C_{k_s},$ then $\alpha_{s}(C)<1/4$ and any  strictly stationary ergodic and aperiodic Markov chain generated by $C=a_1C_1+\cdots+a_k C_k\quad \text{for } \quad 0<a_1,\dots,a_k<1 \quad \text{is}\quad \alpha-\text{mixing}.$ 
\end{theorem}

\subsection{On $\rho$-mixing coefficient}
Bradley (2005) states the following in his Theorem 3.3:
\begin{theorem}\label{Brad2}
Suppose $(X_n,n\in Z)$ is  a  (not  necessarily  stationary) Markov chain. Then each of the following statements holds:

  If $\rho_n<1$ for some $n\ge1$, then $\rho_n\to 0$ at least exponentially fast as $n\to \infty$.

\end{theorem}
Theorem \ref{Brad2} requires $\rho_n<1$, therefore the following is easily established.

\begin{theorem}\label{Theop}
For any set of copulas $C_1(u,v)\cdots C_k (u,v)$, if there exists a subset of copulas $C_{k_1}\cdots C_{k_s}, $ $s\leq k\in \mathbb{N}$ such that $ \rho(\hat{C})<1 \quad \text{for}\quad \hat{C}=C_{k_1}*\cdots*C_{k_s},$ then $\rho_{s}(C)<1$ and any Markov chain generated by $C=a_1C_1+\cdots+a_k C_k\quad \text{for } \quad 0<a_1,\dots,a_k<1 \quad \text{is exponential}\quad \rho-\text{mixing}.$ 
\end{theorem}

\section{Proofs, conclusions and remarks}
\subsection{Proof of Lemma \ref{bound}}
Without loss of generality, we can assume that $(X,Y)$ has cumulative distribution $C(x,y)$. This means that $X$ and $Y$ are uniform on $(0,1)$. Thus, $$P(X_0<x, X_n<y)=C^{n}(x,y)=C^{n-1}*C(x,y).$$ If by induction we assume $c^{n-1}(x,y) \ge K^{n-1}$ on a set of Lebesgue measure 1, then $c^{1}(x,y)=c(x,y)\ge K$ and $C^{n}(x,y)=\int_0^1 C^{n-1}_{,2}(x,t)C_{,1}(t,y)dt$ imply $$c^n(x,y)=\int_0^1c^{n-1}(x,t)c(t,y)dt\ge \int_0^1K^{n-1}Kdt=K^n.$$
Integrating both sides of this inequality gives $C^n(x,y)\ge K^n xy$. 
\subsection{Proof of Theorem \ref{IndXZpart}}
It is enough to conduct this proof for $k=3$ and $s=2$. For random variables $X_1, X_2,X_3$ with distributions $F_1, F_2, F_3$, joint distribution $F$ and copula $C(u,v,w)$, let $\tilde{F}_1, \tilde{F}_2$ be the cumulative distributions of $X_1+Z_1$ and $X_2+Z_2$. Under the assumption that the vectors $X$ and $Z$ are independent, if we denote $u_1=\tilde{F}_1(x_1), u_1=\tilde{F}_2(x_2)$ and $u_3=F_3(x_3)$, then we have by Sklar's  theorem 
$$C_1(u_1, u_2,u_3)= P(X_1+Z_1\leq x_1, X_2+Z_2\leq x_2, X_3\leq x_3).$$
Using independence of $X$ and $Z$ and computing by conditioning, we obtain
$$C_1(u_1, u_2,u_3)=\int_{-\infty}^{\infty}\int_{-\infty}^{\infty}P(X_1\leq x_1-t_1, X_2\leq x_2-t_2, X_3\leq x_3)dG_1(t_1)dG_2(t_2)$$

$$=\int_{-\infty}^{\infty}\int_{-\infty}^{\infty}F(x_1-t_1, x_2-t_2, x_3)dG_1(t_1)dG_2(t_2).$$
The last equality gives the needed fomula, with $x_i=\tilde{F}_i^{-1}(u_1)$ using Sklar's theorem.
The proof of Theorem \ref{XZpart} is similar to that of Theorem \ref{IndXZpart} and the corollaries are special cases.

\subsection{Proof of Theorem \ref{asso}}
Copulas of interest are convex combinations of copulas. The dependence coefficients here are linear functionals of convex combinations of copulas and are equal to $0$ for $\Pi$ and $1$ for $M$, except for $\gamma(C_{\theta,\Pi})=2$ and $\gamma(C_{\theta,M})=3$. So, Theorem \ref{asso} holds.
\subsection{Proof of Theorem \ref{alpha} , Theorem \ref{alpha1} and Theorem \ref{Theop}}
Let $(X_k, k\in \mathbb{N})$ be a copula-based Markov chain generated by $C(u,v)$ and $(\tilde{X}_k, k\in \mathbb{N})$ a Markov chain generated by the perturbation. Recall that $P^n(A\cap B)=P(X_1\in A, X_{n+1}\in B)$ and $\tilde{P}^n(A\cap B)=P(\tilde{X}_1\in A, \tilde{X}_{n+1}\in B)$. When the perturbation uses $\Pi(u,v)$, formula \eqref{Pcop} leads to 

$\tilde{P}^n(A\cap B)-P(A)P(B)=(1-\theta)^n (P^{n}(A\cap B)-P(A)P(B))$. Therefore, $\alpha_n(C_{\theta,\Pi})=(1-\theta)^n \alpha_n(C).$ The rest follows from the fact that $\alpha_n(C)$ is bounded.

When the perturbation uses $M(u,v)$, formula \eqref{Mcop} leads to 

$\displaystyle \tilde{P}^n(A\cap B)-P(A)P(B)=\sum_{i=1}^{n}{n \choose i}\theta^{n-i}(1-\theta)^{i} (P^{n}(A\cap B)-P(A)P(B))+$ 

$+(1-\theta)^n(M(A\cap B)-P(A)P(B))$. 

Let $c$ be the limit of $\alpha_n$ if it exists. Based on the fact that the supremum of a sum is less than the sum of suprema and the limit of $(1-\theta)^n$ is 0, we conclude that 
$\displaystyle |\alpha_n(C_{\theta,M})-c|\leq \sum_{i=1}^{n}{n \choose i}\theta^{n-i}(1-\theta)^{i} |\alpha_i(C)-c|+(1-\theta)^n|\alpha(M)-c|\to \lim_{n\to \infty}|\alpha_n(C)-c|=0$. 

Therefore, when $\lim_{n\to \infty}\alpha_n(C)=0$, we obtain $\lim_{n\to \infty}\alpha_n(C_{\theta,M})=0$.

Theorem \ref{alpha1} is a direct consequence of Theorem \ref{Brad1} and formula \eqref{Pcop}. Here, one of the terms of the sum being strictly less than $1/4$ in absolute value, the entire sum is strictly less than $1/4$. The proof benefits from the fact that we have a Markov chain and the fold product helps derive the needed joint distribution.
As for Theorem \ref{Theop}, the same arguments work with the difference that the correlation is compared to 1 and the idea still works.
\subsection{Conclusion and remarks}
When dealing with Markov chains in the copula framework, the work with mixing coefficients benefits from the fold product of copulas and this helps get more insights as compare to the general case. We have shown the impact of perturbations on the $\alpha$-mixing and $\rho$-mixing structure of Markov chains generated by copulas and will be looking to extend this work to $\psi'$-mixing and beyond.


\begin{thebibliography}{99}
\bibitem{Beare} B.K. Beare (2010). Copulas and Temporal Dependence. {\it Econometrica 78, Issue 1}, 395--410;
\bibitem{Blum} J.R. Blum, D.L. Hanson and L.H. Koopmans (1963). On the strong law of large numbers for a class of stochastic processes. {Z. Wahrscheinlichkeitstheorie und Verw. Gebiete 2, 1--11;}
\bibitem{BR} R.C. Bradley (2007). {\it Introduction to Strong Mixing Conditions. Vol. 1,2}, Kendrick Press;
\bibitem{BR1} R.C. Bradley (2005). Basic Properties of Strong Mixing Conditions. A Survey and Some Open Questions. {\it Probability surveys 2,  107-144;}
\bibitem{BR2} R.C. Bradley (1983). On the $\psi$-mixing condition for stationary random sequences. {\it Transactions of the American Mathematical Society, 276(1) 55--66;}
\bibitem{C} R. Cogburn (1960). Asymptotic properties of stationary sequences. {\it Univ.  Calif. Publ. Statist.} 3 99--146;
\bibitem{DAR} W. F. Darsow, B. Nguyen, E. T. Olsen (1992). {Copulas and Markov processes.} {\it Illinois journal of mathematics  36(4) 600--642};
\bibitem{Gap} V. Gaposhkin (1991). {Moment bounds for integrals of $\rho$-mixing fields}. {\it Theory Probab. Appl. 36, 249--260}; 
\bibitem{I} I.A. Ibragimov (1959). Some limit theorems for stochastic processes stationary in the strict sense. {\it Dokl. Akad. Nauk SSSR} 125 711--714;
\bibitem{IBR} R. Ibragimov (2009). {Copula-based characterizations for higher-order Markov processes.}	{\it Econometric Theory 25 819--846;}
\bibitem{KR} A.N.  Kolmogorov and Yu.A.  Rozanov (1960).  On strong mixing conditions for stationary Gaussian processes. {\it Theor. Probab. Appl.}5  204--208;
\bibitem{KKK}  J. Komornik, M. Komornikova, J. Kalicka (2017). Dependence measures for perturbations of copulas. {\it Fuzzy Sets and Systems} 324 100--116; 
\bibitem{DSF} F. Durante, J.F. Sanchez, M.U. Flores (2013). Bivariate copulas generated by perturbations. {\it Fuzzy Sets and Systems} 228 137--144;
\bibitem{L0}M. Longla, F. Djongreba Ndikwa, M. Muia Nthiani (2021).  Perturbations of copulas and Mixing properties. {\it Journal of the Korean Statistical Society} {doi: 10.1007/s42952-021-00133-5}; 
\bibitem{L1} M. Longla (2015). On mixtures of copulas and mixing coefficients. {\it Journal of Multivariate analysis} 139, 259--265;
\bibitem{L2} M. Longla (2014). On dependence structure of copula-based Markov chains. {\it ESAIM: Probability and Statistics} 18, 570--583;
\bibitem{L3} M. Longla (2013). Remarks on the speed of convergence of mixing coefficients and applications. {\it Statistics and probability letters} 83(10); 2439--2445;
\bibitem{LP} M Longla, M Peligrad (2012). Some aspects of modeling dependence in copula-based Markov chains {\it Journal of Multivariate Analysis} 111, 234--240;
\bibitem{MKK} R. Mesiar, M. Komornikova, J. Komornik (2015). Perturbation of bivariate copula. {\it Fuzzy Sets and Systems} 268 127--140;
\bibitem{Mill}	C. Miller (1994). {Three theorems on $\rho^{*}$-mixing random fields.}  {\it J. Theor. Probab. 7}, 867--882;
\bibitem{Mill1}	C. Miller (1995). { A CLT for the periodograms of a $\rho^{*}$-mixing random field.} {\it Stochastic Process. Appl 60}, 313--330;

\bibitem{N} R.B. Nelsen. {\it An Introduction to Copulas, second edition}, Springer Series in Statistics, Springer-Verlag, New York;

\bibitem{R}  M.  Rosenblatt (1956).  A  central  limit  theorem  and  a  strong  mixing  condition. {\it Proc. Natl. Acad. Sci. USA} 42 43--47;

\bibitem{SK} A. Sklar (1959). Fonctions de répartition \`{a} $n$ dimensions et leurs marges.{\it Publ. Inst. Statist. Univ. Paris}, 8: 229--231;

\end{thebibliography}
\end{document}